\documentclass[11pt,twoside]{article}
\usepackage{amsfonts}
\usepackage{amsmath}
{\markboth{\sl C. Abdelkefi and F. Sassi
}{\sl Besov-Dunkl spaces}}

\newenvironment{pf}{\textit{Proof.}} {\hfill $\Box$}

\begin{document}
\titlepage
\begin{flushleft}
 {\bf{I}}nternational {\bf{J}}ournal of {\bf{P}}ure and {\bf{A}}pplied {\bf{M}}athematics \\
  \small Volume 39, No 4, pp. 475-488, 2007  \;\quad\qquad\qquad\qquad\qquad\qquad\qquad ISSN 1311-8080

\end{flushleft}
\vskip0,5cm
\begin{center}\bf\Large Various characterizations of Besov-Dunkl spaces \end{center}%
\begin{center} {Chokri Abdelkefi$^\S$ and Feriel Sassi $^\ddag$}\footnote{\small Received: April 12, 2007 \hspace{3.3cm} \copyright
\;2007 Academic Publications
 }\end{center}\begin{center} \small
$^\S$$^\ddag$Department of Mathematics, Preparatory Institute of
Engineer Studies of Tunis \\ \small 1089 Monfleury Tunis, Tunisia\\
\small E-mail : chokri.abdelkefi@ipeit.rnu.tn\\\small E-mail :
feriel.sassi@fst.rnu.tn \end{center}
\noindent{\bf Abstract:} In this paper, different characterizations
of the Besov-Dunkl spaces, previously considered in [1, 2, 3, 11],
are given. We provide equivalence between these characterizations,
using the Dunkl translation, the Dunkl transform and the Peetre
$K$-functional.
\vskip1cm
\begin{center}{\bf 1. Introduction}\end{center} \indent $ $
On the real line, we consider the first-order
differential-difference operator defined by
$$\Lambda_\alpha(f)(x) = \frac{df}{dx}(x) + \frac{2\alpha+1}{x}
\left[\frac{f(x) - f(-x)}{2}\right],\quad f \in \mathcal{E}(
\mathbb{R}),\quad \alpha>-\frac{1}{2}\,,$$  which is called  Dunkl
operator. Such operators have been introduced in 1989, by C. Dunkl
in [8]. The Dunkl kernel $E_{\alpha}$ is used to define the Dunkl
transform $\mathcal{F}_{\alpha}$ which was introduced by C. Dunkl in
[9]. R\"osler in [17] shows that the Dunkl kernel verify a product
formula. This allows us to define the Dunkl translation $\tau_{x}$ ,
$x\in\mathbb{R}$. As a result, we have the Dunkl
convolution.\setcounter{page}{475}

There are many ways to define Besov spaces (see [4, 5, 15, 21]).
This paper deals with Besov-Dunkl spaces (see [1, 2, 3, 11]). Let
$1\leq p<+\infty$, $1\leq q \leq + \infty$ and $\beta> 0 $, the
Besov-Dunkl space denoted by $
\mathcal{B}\mathcal{D}_{p,q}^{\beta,\alpha}$ is the subspace of
functions $f \in L^p( \mu_\alpha)$ satisfying
\begin{eqnarray*}
\int^{+ \infty}_0 \left(\frac{w_{p,\alpha}(f,x)}{x^\beta}\right)^q
\, \frac{dx}{x} &<& + \infty \;\;\;\quad if\quad q < +\infty
\end{eqnarray*}  and
\begin{eqnarray*}  \sup_{x\in (0,+\infty)}\frac{w_{p,\alpha}(f,x)}{x^\beta} &<& +\infty\,\quad\quad if \quad q=+\infty ,\end{eqnarray*}
where $ w_{p,\alpha}(f,x)=\displaystyle{\sup_{|t|\leq x}\|\tau_t(f)
- f\|_{p,\alpha}}$ and $\mu_\alpha$ is a weighted Lebesgue measure
on $\mathbb{R}$ (see next section).\\ Put
$\displaystyle{\mathcal{D}_{p,\alpha}}$  the subspace of functions
$f \in L^p( \mu_\alpha)$ such that the distribution function
$\Lambda_{\alpha}f \in L^p( \mu_\alpha).$ $ \mathcal{D}_{p,\alpha }
$ is a Banach space with $\|.\|_{ \mathcal{D}_{p,\alpha} }$ defined
by
\begin{eqnarray*}
\|f\|_{\small{\mathcal{D}_{p,\alpha}}} =\|f\|_{p,\alpha}
+\|\Lambda_{\alpha}f\|_{p,\alpha}.
\end{eqnarray*}
We consider the subspace
$ \mathcal{K} \mathcal{D}_{p,q }^{\beta
,\alpha} $ of functions $f\in L^p( \mu_\alpha)$ satisfying
\begin{eqnarray*}
\int^{+ \infty}_0 \left(\frac{K_{p,\alpha}(f,x)}{x^\beta}\right)^q
\, \frac{dx}{x} &<& + \infty \;\;\;\quad if\quad q < +\infty
\end{eqnarray*}  and
\begin{eqnarray*}  \sup_{x\in (0,+\infty)}\frac{K_{p,\alpha}(f,x)}{x^\beta} &<& +\infty\,\quad\quad if \quad q=+\infty ,\end{eqnarray*}
where $K_{p,\alpha}$ is the Peetre $K$-functional (see[12]) given by
\begin{eqnarray*} K_{p,\alpha}(f,x)=\inf\Big\{\|f_0\|_{p,\alpha}
+x\|\Lambda_{\alpha}f_1\|_{p,\alpha}\,;\, f_0 \in L^p(
\mu_\alpha),\,f_1 \in
\mathcal{D}_{p,\alpha},\;f=f_0+f_1\Big\}.\end{eqnarray*} We denote
by $
  \mathcal{ED}_{p,q }^{\beta ,\alpha}$ the subspace of
functions $f \in L^p( \mu_\alpha)$ satisfying
\begin{eqnarray*}\int^{+ \infty}_1 \left( x^\beta \mathbf{E}_{p,\alpha}(f,x)  \right)^q
\, \frac{dx}{x} &<& + \infty \;\;\;\quad if\quad q < +\infty
\end{eqnarray*}  and
\begin{eqnarray*}  \sup_{x\in (1,+\infty)}x^\beta \mathbf{E}_{p,\alpha}(f,x)  &<& +\infty\,\quad\quad if \quad q=+\infty ,\end{eqnarray*}
 where  $ \mathbf{E}_{p,\alpha}(f,x) =\inf \Big\{\|f-g\|_{p,\alpha }\,;\,
  \; \mbox{supp}\,(\mathcal{F}_\alpha (g))\subset [-x,\;x]\Big\}\,,\;
x>0 .$

Our objective will be to prove that $
\mathcal{B}\mathcal{D}_{p,q}^{\beta,\alpha}=
\mathcal{K}\mathcal{D}_{p,q }^{\beta ,\alpha}$ and when $1\leq p\leq
2,$   $1\leq q <+\infty$, $0<\beta<1$ then $
\mathcal{B}\mathcal{D}_{p,q}^{\beta,\alpha}=
 \mathcal{ED}_{p,q }^{\beta ,\alpha}$.
 \\Analogous results have been
obtained by Betancor, M\'endez and Rodr\'iguez-Mesa in [6] for the Bessel operator on $(0,+\infty)$.

The contents of this paper are as follows. \\In section 2, we
collect some basic definitions and results about harmonic analysis
associated with Dunkl operators .\\ In section 3, we prove the
results about inclusion and coincidence between the spaces
$\mathcal{B}\mathcal{D}_{p,q}^{\beta,\alpha}$ , $ \mathcal{K}\mathcal{D}_{p,q
}^{\beta ,\alpha} $ and $\mathcal{ED}_{p,q }^{\beta ,\alpha}$.

In the sequel $c$ represents a suitable positive constant which is
not necessarily the same in each occurence. Furthermore, we denote
by
 \begin{itemize}
  \item $\mathcal{D}_{\ast}(\mathbb{R})$ the space of even $C^{\infty}$-functions on $\mathbb{R}$
 with compact support.
  \item $\mathcal{S}_{\ast}(\mathbb{R})$ the space of even Schwartz functions
 on $\mathbb{R}$.
 \end{itemize}
\begin{center}{\bf 2. Preliminaries}\end{center}
$ $ Let $\mu_\alpha$ the weighted Lebesgue measure on $\mathbb{R}$
given by
$$d\mu_\alpha(x) = \frac{|x|^{2\alpha+1}}{2^{\alpha +1}\Gamma(\alpha
+1)}dx.$$
For every $1 \leq p \leq + \infty$, we denote by
$L^p(\mu_\alpha)$ the space $L^{p}(\mathbb{R}, d\mu_\alpha)$ and we
use $\|\ \;\|_{p,\alpha}$ as a shorthand for $\|\
\;\|_{L^p(\mu_\alpha)}.$\\
The Dunkl transform $\mathcal{F}_\alpha$
which was introduced by C. Dunkl in [9], is defined for $f \in
L^1(\mu_\alpha)$ by
$$\mathcal{F}_\alpha(f)(x) = \int_{\mathbb{R}}E_\alpha(-ixy)f(y)
d\mu_\alpha(y), \quad x \in \mathbb{R},$$ where for
$\lambda\in\mathbb{C}$, the Dunkl kernel $E_\alpha(\lambda .)$ is
given by $$E_\alpha(\lambda x) = j_\alpha(i\lambda x) +
\frac{\lambda x} {2(\alpha+1)} j_{\alpha+1} (i\lambda x),\quad   x
\in \mathbb{R}, $$ with $j_\alpha$ the normalized Bessel function of
the first kind and order $\alpha$ (see [22]).\\The Dunkl kernel
$E_\alpha(\lambda\; .)$ is the unique solution on $\mathbb{R}$ of
initial problem for the Dunkl operator (see [8]). We have for all
$x,\,y \in \mathbb{R}$,
\begin{eqnarray}|E_\alpha(-ixy)| \leq
1. \end{eqnarray} According to [7], we have
the following results :
 \begin{itemize}
 \item[i)] For all $f \in L^1(\mu_\alpha)$, we have $\|\mathcal{F}_\alpha(f)\|_{\infty,\alpha}
 \leq \|f\|_{1,\alpha}.$\hfill
 \item[ii)] For all $f \in L^1(\mu_\alpha)$ such that $\mathcal{F}_\alpha(f) \in
 L^1(\mu_\alpha)$, we have the inversion formula
 \begin{eqnarray}f(x) = \int_{ \mathbb{R}} E_\alpha (i\lambda x) \mathcal{F}_\alpha(f)(\lambda)d\mu_\alpha
 (\lambda),\; a.e\; x \in \mathbb{R}.\end{eqnarray}
 \item[iii)] For every $f \in L^2(\mu_\alpha)$, we have the Plancherel
 formula
 $$\|\mathcal{F}_\alpha(f)\|_{2,\alpha} = \|f\|_{2,\alpha}.$$
 \end{itemize} For all $x, y, z \in \mathbb{R}$, consider
 \begin{eqnarray}W_\alpha(x,y,z) = \frac{(\Gamma(\alpha+1)^2)}{2^{\alpha-1}
 \sqrt{\pi}\Gamma(\alpha + \frac{1}{2})}
  (1 - b_{x,y,z} + b_{z,x,y} +
 b_{z,y,x}) \Delta_\alpha(x,y,z)  \end{eqnarray}
where
$$b_{x,y,z} = \left\{ \begin{array}{ll}
\frac{x^2 + y^2 -z^2}{2xy} &\mbox{ if } x, y \in \mathbb{R}
\backslash \{0\},\; z \in \mathbb{R}\\
0 &\mbox{ otherwise }
\end{array}\right.$$
and
$$\Delta_\alpha(x,y,z) = \left\{ \begin{array}{ll}
\frac{([(|x| + |y|)^2 - z^2][z^2-(|x| - |y|)^2])^{\alpha -
\frac{1}{2}}}{|xyz|^{2\alpha}} &\mbox{ if } |z|\in S_{x,y}\\
0 &\mbox{ otherwise }
\end{array}\right.$$
where $$S_{x,y} = \Big[||x| - |y||\;,\; |x| + |y|\Big]. $$ The
kernel $W_\alpha $ (see [17]), is even and we have
$$W_\alpha(x,,y,z) = W_\alpha(y,x,z) = W_\alpha(-x,z,y) =
W_\alpha(-z, y, -x)$$ and
$$\int_{\mathbb{R}}|W_\alpha(x,y,z)|d\mu_\alpha(z) \leq
4. $$ In the sequel we consider the signed measure $\gamma_{x,y}$,
on $\mathbb{R}$, given by
\begin{eqnarray} d\gamma_{x,y}(z) = \left\{
\begin{array}{ll} W_\alpha(x,y,z)d\mu_\alpha(z) &\mbox{ if } x, y
\in \mathbb{R} \backslash \{0\}\\ d\delta_x(z) &\mbox{ if } y = 0\\
d\delta_y(z) &\mbox{ if } x = 0.
\end{array}\right.\end{eqnarray}
  For $x, y \in
\mathbb{R}$ and $f$ a continuous function on $\mathbb{R}$, the Dunkl
translation operator $\tau_x$ is given by $$\tau_x(f)(y) =
\int_{\mathbb{R}} f(z) d\gamma_{x,y}(z) .$$ It was shown in [13] that for
$x\in\mathbb{R},$ $\tau_x$ is a continuous linear operator from
$\mathcal{E}( \mathbb{R})$ into itself and for all $f \in
\mathcal{E}(\mathbb{R})$, we have \begin{eqnarray*} \tau_0(f)(x)
=f(x)\,,\;\;\tau_x \circ \tau_y  =  \tau_y
\circ\tau_x\end{eqnarray*}   \begin{eqnarray} \tau_x(f)(y) =
\tau_y(f)(x),\quad \Lambda_\alpha \circ\tau_x=\tau_x\circ
\Lambda_\alpha\,,\;\;\;x,\;y\in \mathbb{R},\end{eqnarray} where
$\mathcal{E}(\mathbb{R})$
denotes the space of $C^\infty$-functions on $\mathbb{R}$.\\
According to [19], the operator $\tau_x$ can be extended to
$L^p(\mu_\alpha) $, $ 1\leq p \leq+\infty$ and for $f \in
L^p(\mu_\alpha)$ we have
\begin{eqnarray} \|\tau_x(f)\|_{p,\alpha} \leq 4
\|f\|_{p,\alpha},\end{eqnarray} and for all $x, \lambda \in
\mathbb{R}$, $f \in L^1(\mu_\alpha)$, we have \begin{eqnarray}
\mathcal{F}_\alpha(\tau_x(f))(\lambda) = E_\alpha(i\lambda x)
\mathcal{F}_\alpha(f)(\lambda).\end{eqnarray} Using the change of
variable $z= (x,y)_\theta=\sqrt{x^2+y^2-2xy \cos\theta}$, we have
also
 \begin{eqnarray} \tau_x(f)(y) =
\int_{0}^\pi\Big[f_e((x,y)_\theta)+\frac{x+y}{(x,y)_\theta}f_o((x,y)_\theta)\Big]
d\nu_\alpha(\theta)\end{eqnarray} where
$$f_e((x,y)_\theta)=f((x,y)_\theta)+f(-(x,y)_\theta)\;,\quad
 f_o((x,y)_\theta)=f((x,y)_\theta)-f(-(x,y)_\theta)$$ and
$$d\nu_\alpha(\theta)=\frac{\Gamma(\alpha+1)}{2\sqrt{\pi}\Gamma(\alpha+\frac{1}{2})}
(1-\cos\theta)\sin^{2\alpha}\theta d\theta .$$
The Dunkl convolution $f\, \ast_\alpha g$ , of two continuous
functions $f$ and $g$ on $\mathbb{R}$ with compact support, is
defined by
$$(f\,\ast_\alpha\, g)(x) = \int_{\mathbb{R}} \tau_x(f)(-y) g(y)
d\mu_\alpha(y),\quad x \in \mathbb{R}.  $$ The convolution
$\ast_\alpha$ is associative and commutative (see [17]).
\\We have the following results (see [18]).
\begin{itemize}\item[i)] Assume
that $p,q, r \in [1, + \infty[$ satisfying $\frac{1}{p} +
\frac{1}{q} = 1 + \frac{1}{r}$ (the Young condition). Then the map
$(f, g) \rightarrow f\,\ast_\alpha \,g$ defined on $C_c(\mathbb{R})
\times C_c( \mathbb{R})$, extends to a continuous map from
$L^p(\mu_\alpha) \times L^q(\mu_\alpha)$ to $L^r(\mu_\alpha)$ and we
have \begin{eqnarray} \|f\,\ast_\alpha \,g\|_{r,\alpha} \leq 4
\|f\|_{p,\alpha} \|g\|_{q,\alpha}. \end{eqnarray}
 \item[ii)] For all $f
\in L^1(\mu_\alpha)$ and $g \in L^2(\mu_\alpha)$, we have
\begin{eqnarray}  \mathcal{F}_\alpha(f\,\ast_\alpha g) =
\mathcal{F}_\alpha(f) \mathcal{F}_\alpha(g)\end{eqnarray}  and for
$f \in L^1(\mu_\alpha)$, $g \in L^p(\mu_\alpha)$ and $1 \leq p <
\infty$, we get \begin{eqnarray} \tau_t (f\,\ast_\alpha\,g) =
\tau_t(f)\,\ast_\alpha g = f\,\ast_\alpha \tau_t(g),\quad t \in
\mathbb{R}.\end{eqnarray}
 \end{itemize}

 \begin{center} {\bf 3. Characterizations of the Besov-Dunkl spaces}\end{center}
    $ $
In this section, we provide equivalence between different characterizations of the
Besov-Dunkl spaces.\\\\
 {\bf Theorem 1.}\label{alp}
{\it Let $1\leq p<+\infty$, $1\leq q \leq +\infty$ and $\beta>0$,
then
 \begin{eqnarray*}\mathcal{B}\mathcal{D}_{p,q}^{\beta,\alpha}= \mathcal{K} \mathcal{D}_{p,q
}^{\beta ,\alpha}. \end{eqnarray*}}
\begin{pf} For $x>0$ and $0<|z|\leq x$, put $\displaystyle{\Theta(x,z)=\frac{1}{2x^{2\alpha+1}}+
\frac{sgn(z)}{2|z|^{2\alpha+1}}}\;.$\\ We start with the proof of the inclusion
 $\mathcal{B}\mathcal{D}_{p,q}^{\beta,\alpha}\subset \mathcal{K} \mathcal{D}_{p,q
}^{\beta ,\alpha} $. For $f\in\mathcal{B}\mathcal{D}_{p,q}^{\beta,\alpha}$ and $x>0$, we take
\begin{eqnarray*} f_1=\frac{1}{x}\int_{-x}^x
\Theta(x,z)\; \tau_z(f)\;d\mu_\alpha(z).
\end{eqnarray*}    Using the Minkowski's
inequality for integrals and (6), we have
\begin{eqnarray*}
\|f_1\|_{p,\alpha}& \leq & \frac{1}{x}\int_{-x}^x    |\Theta(x,z)
|\;\|\tau_z(f)\|_{p,\alpha}\;d\mu_\alpha(z)
\\
&\leq &c\; \frac{\|f\|_{p,\alpha}}{x}\int_{-x}^x   |\Theta(x,z)
|\;d\mu_\alpha(z)\; \leq c\|f\|_{p,\alpha}.
\end{eqnarray*}
By (5) and the generalized Taylor formula with integral remainder
(see[14], Theorem 2, p. 349), we get
\begin{eqnarray*}
\Lambda_\alpha f_1&=&\frac{1}{x}\int_{-x}^x   \Theta(x,z)\;
\tau_z(\Lambda_\alpha
f)\;d\mu_\alpha(z)\\
&=&\frac{1}{x}(\tau_x (f)-f),
\end{eqnarray*}
then we obtain,
\begin{eqnarray}
x\|\Lambda_\alpha f_1\|_{p,\alpha }\leq c \; w_{p,\alpha }(f,x).
\end{eqnarray}
On the other hand, put $f_0=f-2^{\alpha +2}\Gamma(\alpha +2)f_1$, we
can write
\begin{eqnarray*}
f_0 =-\frac{2^{\alpha +2}\Gamma(\alpha +2)}{x}\int_{-x}^x
\Theta(x,z) (\tau_z(f)-f)d\mu_\alpha(z),
\end{eqnarray*}
by the Minkowski's inequality for integrals, we get
\begin{eqnarray}
\|f_0\|_{p,\alpha}& \leq & \frac{c}{x}\int_{-x}^x
|\Theta(x,z)|\;\|\tau_z(f)-f\|_{p,\alpha}\;d\mu_\alpha(z)\nonumber\\
&\leq & c\; \frac{w_{p,\alpha}(f,x)}{x}\int_{-x}^x   |\Theta(x,z)
|\;d\mu_\alpha(z)\nonumber\\ & \leq & c\;w_{p,\alpha}(f,x).
\end{eqnarray}
Hence by (12) and (13), we deduce that
\begin{eqnarray}
K_{p,\alpha }(f,x)\leq c\;w_{p,\alpha}(f,x).
\end{eqnarray}Let prove now the inclusion
$\mathcal{K} \mathcal{D}_{p,q}^{\beta ,\alpha}\subset
\mathcal{B}\mathcal{D}_{p,q}^{\beta,\alpha}$. For $f\in \mathcal{K}
\mathcal{D}_{p,q}^{\beta ,\alpha}$, $x>0$ and $f_0 \in L^p(
\mu_\alpha),\;f_1 \in  \mathcal{D}_{p,\alpha }$ such that
$f=f_0+f_1$, we have by (6)
\begin{eqnarray}w_{p,\alpha}(f_0,x)\leq c\;\|f_0\|_{p,\alpha}\;,
\end{eqnarray}on the other hand, using ([14], Theorem 2) we can write for $t$ such that $|t|\leq x$
\begin{eqnarray*}\tau_t (f_1)-f_1=\int_{-|t|}^{|t|}   \Theta(t,z)\;
\tau_z(\Lambda_\alpha
f_1)\;d\mu_\alpha(z),
\end{eqnarray*}by the Minkowski's inequality for integrals and (6) again, we get
\begin{eqnarray*}
\|\tau_t (f_1)-f_1\|_{p,\alpha}&\leq &\int_{-|t|}^{|t|}   |\Theta(t,z)|\;
\|\tau_z(\Lambda_\alpha
f_1)\|_{p,\alpha}\;d\mu_\alpha(z)\\&\leq& c\; \| \Lambda_\alpha
f_1 \|_{p,\alpha}\int_{-|t|}^{|t|}   |\Theta(t,z)|
\;d\mu_\alpha(z)\\&\leq& c\;|t|\;\| \Lambda_\alpha
f_1 \|_{p,\alpha} \leq c\;x\;\| \Lambda_\alpha
f_1 \|_{p,\alpha}\;,
\end{eqnarray*}then we obtain,
\begin{eqnarray} w_{p,\alpha}(f_1,x)  \leq c\;x\;\| \Lambda_\alpha
f_1 \|_{p,\alpha},
\end{eqnarray}since
\begin{eqnarray*}w_{p,\alpha}(f,x)&\leq& w_{p,\alpha}(f_0,x)+w_{p,\alpha}(f_1,x),
\end{eqnarray*}
by (15) and (16), we deduce that
\begin{eqnarray}w_{p,\alpha}(f,x)&\leq& c\;K_{p,\alpha }(f,x).
\end{eqnarray}Our theorem is proved.
\end{pf}\\\\
 {\bf Theorem 2.}
{\it Let $1\leq p \leq 2,\,1\leq q \leq +\infty$ and $\beta>0$, then
$$ \mathcal{B}\mathcal{D}_{p,q}^{\beta,\alpha}\subset
\mathcal{ED}_{p,q }^{\beta ,\alpha}.$$}
\begin{pf}
Let $f \in\mathcal{B}\mathcal{D}_{p,q}^{\beta,\alpha}$ and $\lambda
,\; x>0$, by (14) and (17) we have
\begin{eqnarray}
w_{p,\alpha}(f,\lambda x)&\leq & c\; K_{p,\alpha }(f,\lambda x)\nonumber \\
&\leq & c \; \max\{1,\lambda \}K_{p,\alpha }(f, x)\nonumber \\
&\leq & c \; \max\{1,\lambda \}w_{p,\alpha }(f, x).
\end{eqnarray}
Choose $\varphi \in \mathcal{S}_{\ast}(\mathbb{R})$ with $ \mbox{supp}\,(\mathcal{F}_\alpha (\varphi )) \subset [-1,1]$ and
 $\displaystyle{\int_\mathbb{R}\varphi (x)d\mu_\alpha(x) =1}.$\\
From (10), we get for $t>0$ $$\mathcal{F}_\alpha (f \ast_\alpha
\varphi_{\frac{1}{t} })=\mathcal{F}_\alpha (f)\mathcal{F}_\alpha
(\varphi_{\frac{1}{t} })$$ where $\varphi_{\frac{1}{t} }(x)=
t^{2(\alpha +1)}\varphi (tx)$, which implies
$\mbox{supp}\,(\mathcal{F}_\alpha (f \ast_\alpha
\varphi_{\frac{1}{t} }))\subset [-t,t]$ and
\begin{eqnarray}
\mathbf{E}_{p,\alpha }(f,t)  \leq \|f-f \ast_\alpha
\varphi_{\frac{1}{t} }\|_{p,\alpha}.
\end{eqnarray}
On the other hand, by the Minkowski's inequality for integrals
\begin{eqnarray*}
\|f-f \ast_\alpha \varphi_{\frac{1}{t} }\|_{p,\alpha} & =
&\Big(\int_\mathbb{R}\Big| f(y)-\int_\mathbb{R}\varphi_{\frac{1}{t}
}(z)\tau_y(f)(z)d\mu_\alpha(z) \Big|^p d\mu_\alpha(y)
\Big)^{1/p}\nonumber\\
& = &\Big(\int_\mathbb{R}\Big| \int_\mathbb{R}\varphi_{\frac{1}{t}
}(z)[f(y)-\tau_z(f)(y)]d\mu_\alpha(z) \Big|^p d\mu_\alpha(y)\Big)^{1/p}\nonumber\\
&\leq &\int_\mathbb{R}|\varphi_{\frac{1}{t}
}(z)|\;\|\tau_z(f)-f\|_{p,\alpha}d\mu_\alpha(z)\nonumber\\
&\leq &\int_\mathbb{R}|\varphi_{\frac{1}{t}
}(z)|\;w_{p,\alpha}(f,|z|)\;d\mu_\alpha(z),
\end{eqnarray*}
using (18), we obtain
\begin{eqnarray}
\|f-f \ast_\alpha \varphi_{\frac{1}{t} }\|_{p,\alpha} & \leq &
c\,w_{p,\alpha}(f,\frac{1}{t})\int_\mathbb{R}|\varphi_{\frac{1}{t}
}(z)|\,(1+t|z|\,)d\mu_\alpha(z)\nonumber\\
& \leq &c\,w_{p,\alpha}(f,\frac{1}{t})\int_\mathbb{R}|\varphi
(z)|\,(1+|z|\,)d\mu_\alpha(z)\nonumber\\ & \leq &
c\,w_{p,\alpha}(f,\frac{1}{t}).
\end{eqnarray}
Thus, (19) and (20) imply
\begin{eqnarray*}\int^{+ \infty}_1 \left( t^\beta
\mathbf{E}_{p,\alpha}(f,t)  \right)^q \, \frac{dt}{t}& \leq &
c\int^{+
\infty}_0 (t^\beta w_{p,\alpha}(f,\frac{1}{t}))^q \, \frac{dt}{t}\nonumber\\
&\leq
 & c\int^{+
\infty}_0\left(\frac{w_{p,\alpha}(f,t)}{t^\beta}\right)^q \,
\frac{dt}{t}\,,\hspace{1cm}\mbox{ if }q<+\infty
\end{eqnarray*}
and the same is true for $q=+\infty$. \\
This completes the proof of
the inclusion.
\end{pf} \\ Now, in order to establish that  $
\mathcal{B}\mathcal{D}_{p,q}^{\beta,\alpha}=
 \mathcal{ED}_{p,q }^{\beta ,\alpha}$ for $1\leq p \leq 2$, $1\leq q<+\infty$ and $0<\beta<1$, we need to show some
 useful results.\\ In the following lemma, we prove a Bernstein-type
 inequality for the Dunkl translation operators. An analogous result has been proved by [6, 10] for the generalized translation operators associated with the
 Bessel operator.\\\\
 {\bf Lemma 1.}
 {\it  For $1\leq p<+\infty$, there exists a constant $c>0$ such that for $h\in L^p( \mu_\alpha)$ an even
 differentiable function on $\mathbb{R}$ with $h'\in L^p(
 \mu_\alpha)$ and $y_1,\; y_2 >0$, we have
 $$\|\tau_{y_1}(h)- \tau_{y_2}(h)\|_{p,\alpha}\leq c\,|y_1-y_2|\;
 \|h'\|_{p,\alpha}.$$}
\begin{pf} Using (8) and the fact that $h$ is even, we can assert
that\\\\
$\|\tau_{y_1}(h)- \tau_{y_2}(h)\|_{p,\alpha}^p$
\begin{eqnarray*}
 &=&\int_\mathbb{R}\Big|[\tau_{y_1}(h)-\tau_{y_2}(h)](x)\Big|^pd\mu_\alpha(x)\\
&=&\int_\mathbb{R}\Big|\int_0^\pi[2h((x,y_1)_\theta)-2h((x,y_2)_\theta)]\;
d\nu_\alpha(\theta)\Big|^pd\mu_\alpha(x) \\
 &\leq
&c\,\int_\mathbb{R}\Big(\int_0^\pi\Big|h((x,y_1)_\theta)-h((x,y_2)_\theta)\Big|^p\;
d\nu_\alpha(\theta)\Big)d\mu_\alpha(x)\\
&\leq&c\,\int_\mathbb{R}\Big(\int_0^\pi\Big|\int^1_0\frac{d}{ds}[h((x,
y_2+s(y_1-y_2))_\theta)]\; ds\Big|^pd\nu_\alpha(\theta)\Big)d\mu_\alpha(x),\\
\end{eqnarray*}since \, $\displaystyle{\frac{d}{ds}\Big|(x, y_2+s(y_1-y_2))_\theta\Big|\leq
|y_1-y_2|\,},$ then we can write\\\\
$\|\tau_{y_1}(h)- \tau_{y_2}(h)\|_{p,\alpha}^p$
\begin{eqnarray*}
 &\leq
&c\,|y_1-y_2|^p\int_\mathbb{R}\int_0^\pi\Big|\int^1_0h'((x,
y_2+s(y_1-y_2))_\theta)\;
ds\Big|^pd\nu_\alpha(\theta)d\mu_\alpha(x),\\
&\leq&c\,|y_1-y_2|^p\int^1_0\int_\mathbb{R}\Big(\int_0^\pi\Big|h'((x,
y_2+s(y_1-y_2))_\theta)\Big|^p\sin^{2\alpha}\theta \; d\theta\Big)
d\mu_\alpha(x)ds.
\end{eqnarray*}
$\displaystyle{\int_\mathbb{R}\Big(\int_0^\pi\Big|h'((x,
y_2+s(y_1-y_2))_\theta)\Big|^p\sin^{2\alpha}\theta \; d\theta\Big)
d\mu_\alpha(x)}$
\begin{eqnarray}
=\int_0^{+\infty} \Big(\int_0^\pi\Big|h'((x,
y_2+s(y_1-y_2))_\theta)\Big|^p\sin^{2\alpha}\theta \; d\theta\Big)
d\mu_\alpha(x)\nonumber\\+\int_{-\infty}^0\Big(\int_0^\pi\Big|h'((x,
y_2+s(y_1-y_2))_\theta)\Big|^p\sin^{2\alpha}\theta \; d\theta\Big)
d\mu_\alpha(x).
\end{eqnarray}
By [20], we have for $x\geq 0$,
\begin{eqnarray} \int_0^\pi\Big|h'((x,
y_2+s(y_1-y_2))_\theta)\Big|^p\sin^{2\alpha}\theta \;
d\theta=c_\alpha\;T_{y_2+s(y_1-y_2)}(|h'|^p)(x)
\end{eqnarray}
where $T_y$, $y\geq 0$ is the generalized translation operator
associated with the Bessel operator and
$c_\alpha=\sqrt{\pi}\,\frac{\Gamma(\alpha+\frac{1}{2})}{\Gamma(\alpha+1)}$.\\
On the other hand, by the change of variable $\theta'=\pi-\theta, $
we get for $x\leq 0$, $\displaystyle{\int_0^\pi\Big|h'((x,
y_2+s(y_1-y_2))_\theta)\Big|^p\sin^{2\alpha}\theta \; d\theta}$
\begin{eqnarray}
&=&\int_0^\pi\Big|h'((-x,
y_2+s(y_1-y_2))_{\theta'})\Big|^p\sin^{2\alpha}\theta'\,
d\theta'\nonumber\\&=&c_\alpha \;T_{y_2+s(y_1-y_2)}(|h'|^p)(-x).
\end{eqnarray}
Then from (21), (22) and (23), we obtain\\
$\displaystyle{\int_\mathbb{R}\Big(\int_0^\pi\Big|h'((x,
y_2+s(y_1-y_2))_\theta)\Big|^p\sin^{2\alpha}\theta \; d\theta\Big)
d\mu_\alpha(x)}$
\begin{eqnarray*}
&=&2c_\alpha\int_0^{+\infty}T_{y_2+s(y_1-y_2)}(|h'|^p)(x)d\mu_\alpha(x)\nonumber\\
&\leq& c\,\int_0^{+\infty}|h'|^p(x)d\mu_\alpha(x) \leq
c\,\|h'\|_{p,\alpha}^p.
\end{eqnarray*}
Hence, we deduce
$$\|\tau_{y_1}(h)- \tau_{y_2}(h)\|_{p,\alpha}\leq c\,|y_1-y_2|\; \|h'\|_{p,\alpha},$$ which proves the result.
\end{pf}\\\\
 {\bf Lemma 2.}
{\it For $1\leq p \leq 2$, there exists a constant $c>0$ such that
for
 any $x>0$, any function $g\in L^p( \mu_\alpha)$ with
 $supp\,(\mathcal{F}_\alpha (g))\subset [-x,\,x]$ and $y_1,y_2>0$, we
have
$$\|\tau_{y_1}(g)- \tau_{y_2}(g)\|_{p,\alpha}\leq c\,x\,|y_1-y_2|\;
\|g\|_{p,\alpha}.$$}
\begin{pf}
Let $g\in \mathcal{S}(\mathbb{R})$ with
$\mbox{supp}\,(\mathcal{F}_\alpha (g))\subset [-x,\,x]$. Choose
$\varphi \in \mathcal{D}_\ast(\mathbb{R})$ such that $\varphi(t)= 1
$ if $|t|\leq 1$ and $\varphi (t)= 0 $ if $ |t|\geq 2$. Then by the
inversion formula (2), we have $\varphi =\mathcal{F}_\alpha (h)$ for
some $h\in \mathcal{S}_\ast(\mathbb{R}).$ Put
$h_x(y)=x^{2(\alpha+1)}h(xy)$ for $y\in\mathbb{R}$, then
$\mathcal{F}_\alpha (h_x)(y)=\varphi(\frac{y}{x})=1$ for $|y|\leq
x.$ Note that $\mbox{supp}\,(\mathcal{F}_\alpha (g))\subset
[-x,\,x]$, then using (1), (7) and (10), we can write
$$\mathcal{F}_\alpha (\tau_{y_1}(g)- \tau_{y_2}(g))=\mathcal{F}_\alpha (h_x\ast_\alpha(\tau_{y_1}(g)- \tau_{y_2}(g))),$$
by (2) and (9), we obtain
\begin{eqnarray*}\tau_{y_1}(g)- \tau_{y_2}(g)&=&h_x\ast_\alpha(\tau_{y_1}(g)-
\tau_{y_2}(g))\\&=&(\tau_{y_1}(h_x)- \tau_{y_2}(h_x))\ast_\alpha
g.\end{eqnarray*}The change of variable $t'=xt$ in (3) gives
$$W_\alpha(xy,xz,t')\;x^{2(\alpha+1)}=W_\alpha( y,
z,t),$$ then from (4), we get
$$d\gamma_{xy,xz}(t')=d\gamma_{y,z}(t) \;\mbox{ and } \tau_{y}(h_x)(z)=x^{2(\alpha+1)}\tau_{xy}(h)(xz).$$
Therefore, using the lemma 1 , we have
\begin{eqnarray*}
\|\tau_{y_1}(g)- \tau_{y_2}(g)\|_{p,\alpha }&\leq & 4\,
\|\tau_{y_1}(h_x)- \tau_{y_2}(h_x)\|_{1,\alpha }\; \|g\|_{p,\alpha}\\
&=&4\,\|\tau_{xy_1}(h)- \tau_{xy_2}(h)\|_{1,\alpha}\;\|g\|_{p,\alpha}\\
&\leq & c\,x\,|y_1-y_2|\; \|h'\|_{p,\alpha}\|g\|_{p,\alpha}\\
&\leq &c\,x\,|y_1-y_2|\; \|g\|_{p,\alpha }.
\end{eqnarray*}
Since $ \mathcal{S}(\mathbb{R})$ is a dense subset of $L^p(
\mu_\alpha)$ for $1\leq p< +\infty$ and by (6), we obtain the
result.
\end{pf}\\\\
{\bf Theorem 3.} {\it Let $1\leq p\leq 2,\,1\leq q < +\infty $ and
$0<\beta<1$, then $$ \mathcal{ED}_{p,q }^{\beta ,\alpha}=
\mathcal{B}\mathcal{D}_{p,q}^{\beta,\alpha}.$$}
\begin{pf}
We have only to show that $ \mathcal{ED}_{p,q }^{\beta
,\alpha}\subset \mathcal{B}\mathcal{D}_{p,q}^{\beta,\alpha}.$ Assume
$f\in \mathcal{ED}_{p,q }^{\beta ,\alpha},$ we can consider $f\neq0$
a.e., then we get
\begin{eqnarray*}
\Big(\int^1_0 (t^{-\beta }w_{p,\alpha}(f,t))^q \,
\frac{dt}{t}\Big)^{\frac{1}{q}}&=&
\Big(\sum^{+\infty}_{n=0}\int^{2^{-n} }_{2^{-n-1}}(t^{-\beta
}w_{p,\alpha}(f,t))^q \, \frac{dt}{t}\Big)^{\frac{1}{q}}\\
&\leq & 2^\beta\Big(\sum^{+\infty}_{n=0}(2^{n\beta}w_{p,\alpha}(f,2^{-n}))^q\Big)^{\frac{1}{q}}\\
&=&2^\beta\sum^{+\infty}_{n=0}\lambda_n
2^{n\beta}w_{p,\alpha}(f,2^{-n}),
\end{eqnarray*}
where
$\displaystyle{\lambda_n=\frac{(2^{n\beta}w_{p,\alpha}(f,2^{-n}))^{\frac{q}{q'}}}{\Big(
\displaystyle{\sum^{+\infty}_{n=0}(2^{n\beta}w_{p,\alpha}(f,2^{-n}))^q\Big)^{\frac{1}{q'}}}}}$\hspace{0.5cm}
with $q'$ the conjugate of $q$.\\
By reasoning as in the proof on ([16], Proposition 3.1, p. 88) and
using the lemma 2, we have for $0<\beta<1$,
\begin{eqnarray*}
\Big(\int^1_0 (t^{-\beta }w_{p,\alpha}(f,t))^q \,
\frac{dt}{t}\Big)^{\frac{1}{q}}&\leq &2^\beta c\Big(
\|f\|_p+(\sum^{+\infty}_{m=1}(2^{m\beta}\mathbf{E}_{p,\alpha}(f,2^{m-1}))^q)^{\frac{1}{q}}\Big)
\end{eqnarray*}
Since $\mathbf{E}_{p,\alpha}(f,t)$ is decreasing in $t$ and by (19),
\begin{eqnarray*}
\Big(\sum^{+\infty}_{m=1}(2^{m\beta}\mathbf{E}_{p,\alpha}(f,2^{m-1}))^q\Big)^{\frac{1}{q}}&=&2^\beta
\mathbf{E}_{p,\alpha}(f,1)+\Big(\sum^{+\infty}_{m=2}(2^{m\beta}\mathbf{E}_{p,\alpha}(f,2^{m-1}))^q\Big)^{\frac{1}{q}}\nonumber\\
&\leq& c\,\Big(\|f\|_p +\Big(\int^{+ \infty}_1 \left( t^\beta
\mathbf{E}_{p,\alpha}(f,t)  \right)^q \,
\frac{dt}{t}\Big)^{\frac{1}{q}}\Big).
\end{eqnarray*}
The result of the two inequalities above is
$$\Big(\int^1_0 (t^{-\beta }w_{p,\alpha}(f,t))^q \,
\frac{dt}{t}\Big)^{\frac{1}{q}}\leq  c\Big( \|f\|_p+\Big(\int^{+
\infty}_1 \left( t^\beta \mathbf{E}_{p,\alpha}(f,t) \right)^q \,
\frac{dt}{t}\Big)^{\frac{1}{q}}\Big).$$ On the other hand, we easily
obtain,
\begin{eqnarray*}
\Big(\int_1^{+\infty} (t^{-\beta }w_{p,\alpha}(f,t))^q \,
\frac{dt}{t}\Big)^{\frac{1}{q}}&\leq & c\,
\|f\|_p\Big(\int_1^{+\infty} t^{-\beta
q-1}\,dt\Big)^{\frac{1}{q}}\,\leq c\, \|f\|_p.
\end{eqnarray*}
Hence, we conclude that
$$\Big(\int^{+\infty}_0 (t^{-\beta }w_{p,\alpha}(f,t))^q \,
\frac{dt}{t}\Big)^{\frac{1}{q}}\leq  c\Big( \|f\|_p+\Big(\int^{+
\infty}_1 \left( t^\beta \mathbf{E}_{p,\alpha}(f,t) \right)^q \,
\frac{dt}{t}\Big)^{\frac{1}{q}}\Big).$$ This completes the proof.
\end{pf}\\
\begin{center}{\bf   References}\end{center}
 \begin{enumerate}
\item C. Abdelkefi and M. Sifi, On the uniform convergence of
partial Dunkl integrals in Besov-Dunkl spaces. Fractional Calculus
and Applied Analysis Vol. 9, N. 1 (2006), 43-56.\item C. Abdelkefi
and M.Sifi, Further results of integrability for the Dunkl
transform. Communication in Mathematical Analysis Vol. 2, N.1
(2007), 29-36.\item C. Abdelkefi and M. Sifi, Characterization of
Besov spaces for the Dunkl operator on the real line. Submitted to
Journal of Inequalities in Pure and Applied Mathematics. \item J. L.
Ansorena and O. Blasco, Characterization of weighted Besov spaces,
Math. Nachr. 171 (1995), 5-17.\item O. V. Besov, On a family of
function spaces in connection with embedding and extensions, Trudy
Mat. Inst. Steklov 60 (1961), 42-81. \item J. J. Betancor, J. M.
M\'endez and L. Rodr\'iguez-Mesa, Espacios de Besov Asociados a la
Transformaci\'on de Hankel, Servicio de Publicaciones, Universidad
de la Rioja, Logro\~no, Spain, 2001.\item M. F. E de Jeu, The Dunkl
transform , Inv. Math. 113 (1993) 147-162.\item C.F.
Dunkl, Differential-difference operators associated to reflection\\
groups, Trans.Amer. Math. Soc. 311, No1, (1989), 167-183.\item C.F.
Dunkl, Hankel transforms associated to finite reflection groups  in
: Proc. of special session on hypergeometric functions on domains of
positivity, Jack polynomials and applications. proceedings, Tampa
1991, Contemp. Math. 138 (1992), 123-138.\item J. Gosselin and K.
Stempak, A weak type estimate for Fourier-Bessel multipliers, Proc.
Amer. Math. Soc. 106 (1989), 655-662.\item L. Kamoun, Besov-type
spaces for the Dunkl operator on the real line. Journal of
Computational and Applied Mathematics 199, N. 1 (2007), 56-67.
\item J. L\"ofstr\"om and J. Peetre, Approximation theorems connected with generalized
translations, Math. Ann. 181 (1969), 255-268.   \item M.A. Mourou,
Transmutation operators associated with a Dunkl-type
differential-difference operator on the real line and certain of
their applications, Integral Transforms Spec. Funct. 12, No 1
(2001), 77-88.\item M.A. Mourou, Taylor series associated with a
differential-difference operator on the real line. Journal of
Computational and Applied Mathematics, Vol. 153 (2003), 343-354.
\item J. Peetre, New
thoughts on Besov spaces, Duke Univ. Math. Series, Durham, NC,
1976.\item A. Pelczy\'nski and M. Wojciechowski, Molecular decompositions and embedding
theorems for vector-valued Sobolev spaces with gradient norm, Studia Math. 107 (1993), 61-100.
 \item M. R\"osler, Bessel-Type signed hypergroup on
$\mathbb{R}$, in Probability measure on groups and related
structures, Proc. Conf. Oberwolfach, (1994), H. Heyer and A.
Mukherjea (Eds) World Scientific Publ, 1995, 292-304.
\item F. Soltani, $L^p$-Fourier multipliers for the Dunkl operator
on the real line, J. Funct. Analysis, 209 (2004), 16-35.
\item S. Thangavelyu and Y. Xu, Convolution operator and maximal
function for Dunkl transform. J. Anal. Math., vol. 97, pp. 25-55
(2005).\item K. Trim\`eche, Generalized Harmonic Analysis and
Wavelets Packets. Gordon and Breach Science Publichers 2001
OPA.\item H. Triebel, Theory of function spaces, Monographs in
Math., Vol. 78, Birkh\"auser, Verlag, Basel, 1983.\item G.N. Watson,
A treatise on the theory of Bessel functions, Camb. Univ. Press,
Cambridge 1966.
\end{enumerate}
\end{document}